\documentclass[fleqn]{zzz}
\usepackage{graphics,epsfig,color}
\usepackage{hyperref}

\newcommand{\N}{{\mathbf N}}

\newcommand{\C}{{\mathbf C}}

\def\cprime{$'$}

\makeatletter
\def\eqnarray{\stepcounter{equation}\let\@currentlabel=\theequation
\global\@eqnswtrue
\tabskip\@centering\let\\=\@eqncr
$$\halign to \displaywidth\bgroup\hfil\global\@eqcnt\z@
  $\displaystyle\tabskip\z@{##}$&\global\@eqcnt\@ne
  \hfil$\displaystyle{{}##{}}$\hfil
  &\global\@eqcnt\tw@ $\displaystyle{##}$\hfil
  \tabskip\@centering&\llap{##}\tabskip\z@\cr}

\def\endeqnarray{\@@eqncr\egroup
      \global\advance\c@equation\m@ne$$\global\@ignoretrue}

\def\@yeqncr{\@ifnextchar [{\@xeqncr}{\@xeqncr[5pt]}}
\makeatother

\begin{document}

\renewcommand{\PaperNumber}{***}

\FirstPageHeading

\ShortArticleName{Some dual definite integrals for Bessel functions}

\ArticleName{Some dual definite integrals for Bessel functions}

\Author{Howard S.~COHL,${}^{*}$ Sean J.~NAIR,${}^{\dag}$ and Rebekah M.~PALMER\,${}^{\ddag}$}

\AuthorNameForHeading{H.~S.~Cohl, S.~J.~Nair, \& R.~Palmer}

\Address{${}^{*}$~Applied and Computational Mathematics Division,
National Institute of Standards and Technology,
Gaithersburg, MD, 20899-8910, USA
} 
\EmailD{howard.cohl@nist.gov} 
\URLaddressD{http://www.nist.gov/itl/math/msg/howard-s-cohl.cfm}

\Address{${}^{\dag}$~Mathematics, Science, and Computer Science Magnet Program,
Montgomery Blair High School, Silver Spring, MD, 20901, USA
}
\EmailD{s.j.nair2@gmail.com} 

\Address{$^\ddag$~Department of Mathematics, Johns Hopkins University, 
Baltimore, MD 21218, USA
}
\EmailD{rmaepalmer4@gmail.com} 


\ArticleDates{Received XX July 2012 in final form ????; Published online ????}

\Abstract{
Based on known definite integrals of Bessel functions of the first kind, 
we obtain exact solutions to unknown definite integrals using 
the method of integral transforms from Hankel's transform.} 


\Classification{26A42; 33C05;  33C10; 33C45; 35A08}


\section{Introduction}

In Cohl (2012) \cite{CohlVolkmerDefInt}, orthogonality and Hankel's transform are 
used to generate solutions to new definite integrals based on known integrals.
In this paper, we use the method of integral transforms 
to create new integrals from a variety of known integrals containing Bessel
functions of the first kind $J_\nu$. In this method, we use the closure relation for
Bessel functions of the first kind to generate a guess for a function to use
in Hankel's transform. This guess may be incorrect if it does not satisfy the first condition
for Hankel's transform, in which case a new definite integral is not generated. If 
the guess satisfies the condition, restrictions on $\nu$ must then be adjusted
to satisfy the other condition of Hankel's transform.
For the functions used in this paper, superscripts and subscripts refer to lists of
parameters. Any exceptions to this will be clear from context.

As far as we are aware, the 40 definite integrals over Bessel functions of the first kind
that we present in this manuscript, do not currently appear in the literature. An extension 
of the survey presented in this manuscript can be used to mechanically compute new definite
integrals from pre-existing definite integrals over Bessel functions of the first kind.

\subsection{Application of Hankel's transform}
We use the following result where for $x\in(0,\infty)$ we define
\[
F(r\pm 0):=\lim_{x\to r\pm}F(x);
\]
see Watson (1944) \cite[p.~456]{Watson}:

\begin{theorem}
\label{2:Hankel}
Let $F:(0,\infty)\to\C$ be such that
\begin{equation}\label{2:cond}
\int_0^\infty \sqrt{x}\,|F(x)|\,dx<\infty,
\end{equation}
and let $\nu \ge -\frac12$. Then
\begin{equation}\label{2:Hankel2}
 \frac12(F(r+0)+F(r-0))=\int_0^\infty uJ_\nu(ur)\int_0^\infty xF(x)J_\nu(ux)\,dx\,du
\end{equation}
provided that the positive number $r$ lies inside an interval in which $F(x)$ has 
finite variation.  
\end{theorem}

The effort described in this paper was motivated by the large collection of Bessel
function definite integrals which exist in the book ``Table of Integrals, 
Series, and Products'' \cite{Grad}.  Not counting 
Theorem \ref{Heron} (which stands alone), the method of integral transforms was applied 
to the Bessel function definite integrals appearing in Sections 6.51 and 6.52 of
\cite{Grad}. This method can be applied to many definite integrals appearing in the 
remainder of sections appearing in Sections 6.5-6.7 of \cite{Grad}.

For the definite integrals presented in this manuscript, we have directly verified that 
(\ref{2:cond}) is satisfied. This is easily accomplished by analyzing the behavior of 
the integrands in a small neighborhood of the endpoints $\{0,\infty\}$.
For this paper, this technique produced 30 theorems including 40 definite integrals which
are given below.  The method of integral transforms does not always succeed 
in producing new definite integrals because the conditions on the Hankel transform 
(\ref{2:cond}) is not satisfied. Some cases of this are shown in Section \ref{AttemptFailure}.

\section{Polynomial, rational, algebraic, and power functions}

\begin{theorem}
\label{Heron}
Let $b$,$c>0$, $\nu>-\frac12$, $t\in\C\setminus (-\infty,0]$. Then
\begin{eqnarray}
\int_0^\infty
(\Delta(a,b,c))^{2\nu-1}
a^{1-\nu} J_\nu(at) da
=2^{1-\nu} \sqrt{\pi} \; \Gamma\left(\nu+\tfrac{1}{2}\right) 
\left(\tfrac{bc}{t} \right)^\nu J_\nu(bt) J_\nu(ct).
\end{eqnarray}
$\Delta:[0,\infty)^3\to[0,\infty)$ (Heron's formula \cite{Ostermannetal12}), 
defined by 
\[
\Delta(a,b,c):=\sqrt{s(s-a)(s-b)(s-c)},
\] 
where $s=(a+b+c)/2$, is the area of a triangle with sides of length a, b, and c. 
\end{theorem}
\medskip

\noindent {\bf Proof.}\quad We apply Theorem \ref{2:Hankel} to the function
$F_\nu^{b,c}:(0,\infty)\to\C$ defined by
\[
F_\nu^{b,c}(t):=2^{1-\nu} \sqrt{\pi} \; \Gamma\left(\nu+\tfrac{1}{2}\right)
\left(\tfrac{bc}{t} \right)^\nu J_\nu(bt) J_\nu(ct),
\]
where $\Gamma:\C\setminus -\N_0\to\C$ is Euler's gamma function is defined in \cite[(5.2.1)]{NIST},
and $J_\nu:\C\setminus(-\infty,0]\to\C$, (order) $\nu\in\C$, is the Bessel 
function of the first kind defined in \cite[(10.2.2)]{NIST}. 
The desired result is obtained from Sonine's formula \cite{Stempak}
\[
\int_0^\infty J_\nu(at) J_\nu(bt) J_\nu(ct) t^{1-\nu} dt
=
\frac{2^{\nu-1}(\Delta(a,b,c))^{2\nu-1}}{\sqrt{\pi} \;
\Gamma\left(\nu+\tfrac{1}{2}\right) (abc)^\nu},
\]
where ${\rm Re}\,a>0$, $b$,$c>0$, ${\rm Re}\,\nu>-\frac12$.
$\hfill\blacksquare$
\medskip

\begin{theorem}
Let $\nu>\frac12$, $\mu>0$, $\alpha$,$\beta>0$, $z\in\C\setminus(-\infty,0]$. Then
\begin{eqnarray}
\int_{0}^{\alpha} b^{\nu} J_{\nu-1}(bz)\, db
=
  \alpha^{\nu} z^{-1} J_{\nu}(\alpha z),
\end{eqnarray}
\begin{eqnarray}
\int_{\beta}^{\infty} a^{1-\mu} J_{\mu}(az)\, da
=
  \beta^{1-\mu} z^{-1} J_{\mu-1}(\beta z).
\end{eqnarray}
\end{theorem}
\medskip

\noindent {\bf Proof.}\quad By applying Theorem \ref{2:Hankel} to the functions
$F_\nu^\alpha:(0,\infty)\to\C$ and $G_\mu^\beta:(0,\infty)\to\C$ defined by 
$
F_\nu^\alpha(x):=\alpha^{\nu} x^{-1} J_{\nu}(\alpha x)
$,
$
G_\mu^\beta(x):=\beta^{1-\mu} x^{-1} J_{\mu-1}(\beta x)
$,
we obtain the desired results from the known integral \cite[(6.512.3)]{Grad}
\[
\int_{0}^{\infty} J_{\nu}(\alpha x) J_{\nu-1}(\beta x) dx=\left\{
\begin{array}{ll}
    \alpha^{-\nu} \beta^{\nu-1} & $if$\,\, \beta<\alpha, \\[.2cm]
    (2\beta)^{-1} & $if$\,\, \beta=\alpha, \\[.2cm]
    0 & $if$\,\, \beta>\alpha, \\[.2cm]
\end{array}
\right.
\]
where ${\rm Re}\,\nu>0$.
$\hfill\blacksquare$
\medskip

\begin{theorem}
Let $\nu\geq-\frac12$, $z\in\C\setminus(-\infty,0]$. Then
\begin{eqnarray}
\int_0^\infty \frac{c^{\nu+1}}{1+c^2} J_\nu(cz) dc
=K_\nu(z).
\end{eqnarray}
\end{theorem}
\medskip

\noindent {\bf Proof.}\quad We are given the integral \cite[(6.521.2)]{Grad}
\[
\int_{0}^{\infty} x K_{\nu} (ax) J_{\nu} (bx) dx
=
\frac{b^{\nu}}{a^{\nu} (b^2+a^2)},
\]
where  ${\rm Re}\,a>0$, $b>0$, ${\rm Re}\,\nu>-1$.
By applying Theorem \ref{2:Hankel} to the function
$F_\nu^a:(0,\infty)\to\C$ defined by
$
F_\nu^a(x):=a^{\nu} K_{\nu} (ax)
$,
we obtain the following integral
\[
\int_{0}^{\infty} \frac{b^{\nu+1}}{b^2+a^2} J_{\nu} (bx) db
=
  a^{\nu} K_{\nu} (ax),
\]
where ${\rm Re}\,a>0$, $\nu\geq-\frac12$, $x\in\C\setminus(-\infty,0]$.
With the substitutions $z=ax$ and $c=b/a$, we obtain the desired result.
$\hfill\blacksquare$
\medskip

Note that when the method of integral transforms is applied to 
$
F_a(x):=a K_1(ax)
$
given the integral \cite[(6.521.7)]{Grad}
\[
\int_0^\infty x K_1(ax) J_1(bx)
= \frac{b}{a(a^2+b^2)},
\]
where $a>0$, $b>0$,
we obtain the integral generated from \cite[(6.521.2)]{Grad} when $\nu=1$.

\begin{theorem}
Let $z\in\C\setminus(-\infty,0]$. Then
\begin{eqnarray}
\int_0^\infty \frac{c}{(1+c^2)^2} J_0(cz)\, dc
= \frac{z}{2} K_1(z).
\end{eqnarray}
\end{theorem}
\medskip

\noindent {\bf Proof.}\quad We are given the integral \cite[(6.521.12)]{Grad}
\[
\int_{0}^{\infty} x^2 K_{1}(ax) J_{0} (bx)
=
\frac{2a}{(a^2+b^2)^2},
\]
where $a>b>0$.
By applying Theorem \ref{2:Hankel} to the function
$F_a:(0,\infty)\to\C$ defined by
$
F_a(x):=\frac{x}{2a} K_{1}(ax)
$,
we obtain the following integral
\[
\int_{0}^{\infty} \frac{b J_{0} (bx)}{(a^2+b^2)^2} db
=
\frac{x}{2a} K_{1}(ax),
\]
where $a>0$, $x\in\C\setminus(-\infty,0]$.
With the substitutions $z=ax$ and $c=b/a$, we obtain the desired result.
$\hfill\blacksquare$
\medskip

\begin{theorem}
Let $z\in\C\setminus(-\infty,0]$. Then 
\begin{eqnarray}
\int_0^\infty \frac{c^2}{(1+c^2)^2} J_1(cz) dc
=\frac{z}{2}K_0(z).
\end{eqnarray}
\end{theorem}
\medskip

\noindent {\bf Proof.}\quad We are given the integral \cite[(6.521.12)]{Grad}
\[
\int_{0}^{\infty} x^2 K_{0}(ax) J_{1} (bx) dx
=
\frac{2b}{(a^2+b^2)^2},
\]
where $a$,$b>0$.
By applying Theorem \ref{2:Hankel} to the function
$F_a:(0,\infty)\to\C$ defined by
$
F_a(x):=\frac{x}{2} K_{0}(ax)
$,
we obtain the following integral
\[
\int_{0}^{\infty} \frac{b^2 J_{1} (bx)}{(a^2+b^2)^2} db
=
\frac{x}{2} K_{0}(ax),
\]
where $a>0$, $x\in\C\setminus(-\infty,0]$.
With the substitutions $z=ax$ and $c=b/a$, we obtain the desired result.
$\hfill\blacksquare$

\begin{theorem}
Let $\gamma>0$, $\nu\geq-\frac12$,  ${\rm Re}\,\alpha>|{\rm Im}\,\beta|$, $z\in\C\setminus(-\infty,0]$. Then
\begin{eqnarray}
\int_0^\infty b J_\nu(bz) \frac{l_1^\nu}{l_2^\nu(l_2^2-l_1^2)} db
=K_0(\alpha z) J_\nu(\gamma z),
\end{eqnarray}
\begin{eqnarray}
\int_0^\infty c J_\nu(cz) \frac{l_1^\nu}{l_2^\nu(l_2^2-l_1^2)} dc
=K_0(\alpha z) J_\nu(\beta z),
\end{eqnarray}
where $l_1$ and $l_2$ are defined as 
\begin{equation}\label{2:l1}
l_1=\frac{1}{2} \left[\sqrt{(b+c)^2+a^2} - \sqrt{(b-c)^2+a^2}\right],
\end{equation}
\begin{equation}\label{2:l2}
l_2=\frac{1}{2} \left[\sqrt{(b+c)^2+a^2} + \sqrt{(b-c)^2+a^2}\right].
\end{equation}
\end{theorem}
\medskip

\noindent {\bf Proof.}\quad By applying Theorem \ref{2:Hankel} to the function
$F_\nu^{a,c}:(0,\infty)\to\C$ and $G_\nu^{a,b}:(0,\infty)\to\C$ defined by
$
F_\nu^{a,c}(x):= K_0(ax) J_\nu(cx)
$,
$
G_\nu^{a,b}(x):=K_0(ax) J_\nu(bx)
$,
we obtain the desired result from the known integral \cite[(6.522.12)]{Grad}
\[
\int_0^\infty x K_0(ax) J_\nu(bx) J_\nu(cx) dx
= \frac{l_1^\nu}{l_2^\nu(l_2^2-l_1^2)},
\]
where $c>0$,  ${\rm Re}\,\nu>-1$,  ${\rm Re}\,a>|{\rm Im}\,b|$. 
$\hfill\blacksquare$

\begin{theorem}
Let ${\rm Re}\,b > {\rm Re}\,a$, $z\in\C\setminus(-\infty,0]$. Then
\begin{eqnarray}
\int_0^\infty c J_0(cz) (a^4+b^4+c^4-2a^2b^2+2a^2c^2+2b^2c^2)^{-1/2} dc
= I_0(az)K_0(bz).
\end{eqnarray}
Let ${\rm Re}\,b > {\rm Re}\,c$, $z\in\C\setminus(-\infty,0]$. Then
\begin{eqnarray}
\int_0^\infty \frac{a J_0(az)}{l_2^2-l_1^2} da
=I_0(cz)K_0(bz).
\end{eqnarray}
where $l_1$ and $l_2$ are defined in $(\ref{2:l1})$ and $(\ref{2:l2})$ respectively.
\end{theorem}
\medskip

\noindent {\bf Proof.}\quad By applying Theorem \ref{2:Hankel} to the function
$F_a^b:(0,\infty)\to\C$ and $G_c^b:(0,\infty)\to\C$ defined by
$
F_a^b(x):= I_0(ax)K_0(bx)
$,
$
G_c^b(x):=I_0(cx)K_0(bx)
$,
we obtain the desired results from the known integrals (see \cite[(6.522.4)]{Grad})
\[
\int_0^\infty x I_0(ax) K_0(bx) J_0(cx) dx
=(a^4+b^4+c^4-2a^2b^2+2a^2c^2+2b^2c^2)^{-1/2},
\]
where ${\rm Re}\,b > {\rm Re}\,a$, $c>0$, and
\[
\int_0^\infty x I_0(cx) K_0(bx) J_0(ax) dx
=\frac{1}{l_2^2-l_1^2},
\]
where ${\rm Re}\,b > {\rm Re}\,c$, $a>0$,
respectively.
$\hfill\blacksquare$
\medskip

\begin{theorem}
Let $\gamma>0$, $\nu\geq-\frac12$, ${\rm Re}\,\alpha>|{\rm Im}\,\beta|$, $z\in\C\setminus(-\infty,0]$. Then
\begin{eqnarray}
\int_0^\infty \frac{b^{\nu+1}}{(l_2^2-l_1^2)^{2\nu+1}} J_\nu(bz) db
=\frac{z^\nu(\alpha\gamma)^{-\nu} \sqrt{\pi}}{2^{3\nu}\,\Gamma(\nu
+\frac{1}{2})} K_\nu(\alpha z) J_\nu(\gamma z),
\end{eqnarray}
\begin{eqnarray}
\int_0^\infty \frac{c^{\nu+1}}{(l_2^2-l_1^2)^{2\nu+1}} J_\nu(cz) dc
=\frac{z^\nu(\alpha\beta)^{-\nu} \sqrt{\pi}}{2^{3\nu}\,\Gamma(\nu+\frac{1}{2})} K_\nu(\alpha z) J_\nu(\beta z),
\end{eqnarray}
where $l_1$ and $l_2$ are defined in $(\ref{2:l1})$ and $(\ref{2:l2})$, respectively.
\end{theorem}
\medskip

\noindent {\bf Proof.}\quad By applying Theorem \ref{2:Hankel} to the functions
$F_\nu^{a,c}:(0,\infty)\to\C$, $G_\nu^{a,b}:(0,\infty)\to\C$, defined by 
\[
F_\nu^{a,c}(x):=\frac{x^\nu(ac)^{-\nu} \sqrt{\pi}}{2^{3\nu}\,\Gamma(\nu+\frac{1}{2})} K_\nu(ax) J_\nu(cx),
\]
\[
G_\nu^{a,b}(x):=\frac{x^\nu(ab)^{-\nu} \sqrt{\pi}}{2^{3\nu}\,\Gamma(\nu+\frac{1}{2})} K_\nu(ax) J_\nu(bx),
\]
we obtain the desired results from the known integral \cite[(6.522.15)]{Grad}
\[
\int_0^\infty x^{\nu+1} J_\nu(bx) K_\nu(ax) J_\nu(cx) dx
= \frac{2^{3\nu}(abc)^\nu\,\Gamma(\nu+\frac{1}{2})}{\sqrt{\pi}\,(l_2^2-l_1^2)^{2\nu+1}},
\]
where ${\rm Re}\,a>|{\rm Im}\,b|$, $c>0$.
$\hfill\blacksquare$
\medskip

\begin{theorem}
Let $\gamma>0$, ${\rm Re}\,\beta\geq|{\rm Im}\,\alpha|$, ${\rm Re}\,\alpha>0$, 
${\rm Re}\,p>|{\rm Re}\,q|$, ${\rm Re}\,q>0$, $z\in\C\setminus(-\infty,0]$. Then
\begin{eqnarray}
\int_0^\infty 2a^2 J_1(az) (a^2+\beta^2-\gamma^2)\left[(a^2+\beta^2+\gamma^2)^2-4a^2\gamma^2\right]
^{-3/2} da
=z K_0(\beta z) J_0(\gamma z),
\end{eqnarray}
\begin{eqnarray}
\int_0^\infty c J_0(cz) (\alpha^2+\beta^2-c^2)\left[(\alpha^2+\beta^2+c^2)^2-4\alpha^2c^2\right]
^{-3/2} dc
=\frac{z}{2\alpha} J_1(\alpha z) K_0(\beta z).
\end{eqnarray}
\begin{eqnarray}
\int_0^\infty J_1(bz) \frac{2b^2 (p^2+b^2-\gamma^2)}{(l_2^2-l_1^2)^3} db
=z K_0(p z) J_0(\gamma z),
\end{eqnarray}
\begin{eqnarray}
\int_0^\infty J_0(cz) \frac{c (p^2+q^2-c^2)}{(l_2^2-l_1^2)^3} dc
=\frac{z}{2q} J_1(q z) K_0(p z),
\end{eqnarray}
where $l_1$ and $l_2$ are defined in $(\ref{2:l1})$ and $(\ref{2:l2})$, respectively.
\end{theorem}

\noindent {\bf Proof.}\quad By applying Theorem \ref{2:Hankel} to the functions
$F_b^c:(0,\infty)\to\C$, $G_a^b:(0,\infty)\to\C$, $H_a^c:(0,\infty)\to\C$, 
$I_b^a:(0,\infty)\to\C$ defined by 
$
F_b^c(x):=x K_0(bx) J_0(cx)
$,
$
G_a^b(x):=\frac{x}{2a} J_1(ax) K_0(bx)
$,
$
H_a^c(x):=x K_0(ax) J_0(cx)
$,
$
I_b^a(x):=\frac{x}{2b} J_1(bx) K_0(ax)
$,
we obtain the desired results from the known integrals (see \cite[(6.525.1)]{Grad})
\[
\int_0^\infty x^2 J_1(ax) K_0(bx) J_0(cx) dx
=2a(a^2+b^2-c^2)\left[(a^2+b^2+c^2)^2-4a^2c^2\right]^{-3/2},
\]
where $c>0$, ${\rm Re}\,b\geq|{\rm Re}\,a|$, ${\rm Re}\,a>0$, 
\[
\int_0^\infty x^2 J_1(bx) K_0(ax) J_0(cx) dx
= \frac{2b (a^2+b^2-c^2)}{(l_2^2-l_1^2)^3},
\]
where $c>0$, ${\rm Re}\,a>|{\rm Im}\,b|$, ${\rm Re}\,b>0$. 
$\hfill\blacksquare$

\begin{theorem}
Let ${\rm Re}\,a>0$, $\nu\geq-\frac12$, $z\in\C\setminus(-\infty,0]$. Then
\begin{eqnarray}
\int_0^\infty \frac{J_\nu(bz)}{\sqrt{b^2+4a^2}} db
= I_{\nu/2}(az) K_{\nu/2}(az).
\end{eqnarray}
\end{theorem}
\medskip

\noindent {\bf Proof.}\quad By applying Theorem \ref{2:Hankel} to the function
$F_\nu^a:(0,\infty)\to\C$ defined by 
\[
F_\nu^a(x):=I_{\nu/2}(ax) K_{\nu/2}(ax),
\]
we obtain the desired result from the known integral \cite[(6.522.9)]{Grad}
\[
\int_0^\infty x I_{\nu/2}(ax) K_{\nu/2}(ax) J_\nu(bx) dx
= b^{-1}(b^2+4a^2)^{-1/2},
\]
where $b>0$, ${\rm Re}\,a>0$, ${\rm Re}\,\nu>-1$.
$\hfill\blacksquare$

\begin{theorem}
Let $a>0$, $\nu\geq-\frac12$, $z\in\C\setminus(-\infty,0]$. Then
\begin{eqnarray}
\int_{2a}^\infty \frac{J_\nu(bz)}{\sqrt{b^2-4a^2}} db
=-\frac{\pi}{2} J_{\nu/2}(az)Y_{\nu/2}(az).
\end{eqnarray}
\end{theorem}
\medskip

\noindent {\bf Proof.}\quad By applying Theorem \ref{2:Hankel} to the function
$F_\nu^a:(0,\infty)\to\C$ defined by 
\[
F_\nu^a(x):=-\frac{\pi}{2} J_{\nu/2}(ax)Y_{\nu/2}(ax),
\]
we obtain the desired result from the known integral \cite[(6.522.10)]{Grad}
\[
\int_0^\infty x J_{\nu/2}(ax)Y_{\nu/2}(ax) J_\nu(bx) dx
=
\left\{
\begin{array}{ll}
    0 & $if$\,\, 0<b<2a, \\[.2cm]
    -2\pi^{-1}b^{-1}(b^2-4a^2)^{-1/2} & $if$\,\, 2a<b,
\end{array}
\right.
\]
where ${\rm Re}\,\nu>-1$.
$\hfill\blacksquare$

\begin{theorem}
Let ${\rm Re}\,a>0$, $\nu\geq-\frac12$, ${\rm Re}\,\mu\geq\frac32$, $z\in\C\setminus(-\infty,0]$. Then
\begin{eqnarray}
\int_0^\infty \frac{J_\nu(bz)}{\sqrt{b^2+4a^2}} \left[b+(b^2+4a^2)^{1/2}\right]^\mu db
=2^\mu a^\mu I_{(\nu-\mu)/2}(az) K_{(\nu+\mu)/2}(az).
\end{eqnarray}
\end{theorem}
\medskip

\noindent {\bf Proof.}\quad By applying Theorem \ref{2:Hankel} to the function
$F_\nu^{\mu,a}:(0,\infty)\to\C$ defined by 
\[
F_\nu^{\mu,a}(x):=2^\mu a^\mu I_{(\nu-\mu)/2}(ax) K_{(\nu+\mu)/2}(ax),
\]
we obtain the desired result from the known integral \cite[(6.522.12)]{Grad}
\[
\int_0^\infty x I_{(\nu-\mu)/2}(ax) K_{(\nu+\mu)/2}(ax) J_\nu(bx) dx
= 2^{-\mu} a^{-\mu} b^{-1} (b^2+4a^2)^{-1/2} \left[b+(b^2+4a^2)^{1/2}\right]^\mu,
\]
where ${\rm Re}\,a>0$, $b>0$, ${\rm Re}\,\nu>-1$, ${\rm Re}\,(\nu-\mu)>-2$.
$\hfill\blacksquare$
\medskip

\begin{theorem}
Let $|{\rm Re}\,a|<{\rm Re}\,b$, $x\in\C\setminus(-\infty,0]$. Then
\begin{eqnarray}
\int_0^\infty c J_0(cx) (b^2+c^2-a^2) \left[(a^2+b^2+c^2)^2-4a^2b^2\right]
^{-3/2} dc
= \frac{x}{2b} I_0(ax) K_1(bx).
\end{eqnarray}
\end{theorem}
\medskip

\noindent {\bf Proof.}\quad By applying Theorem \ref{2:Hankel} to the function
$F_a^b:(0,\infty)\to\C$ defined by 
$
F_a^b(x):= \frac{x}{2b} I_0(ax) K_1(bx),
$
we obtain the desired result from the known integral \cite[(6.525.2)]{Grad}
\[
\int_0^\infty x^2 I_0(ax) K_1(bx) J_0(cx) dx
= 2b(b^2+c^2-a^2) \left[(a^2+b^2+c^2)^2-4a^2b^2\right]^{-3/2},
\]
where ${\rm Re}\,b>|{\rm Re}\,a|$, $c>0$.
$\hfill\blacksquare$

\section{Bessel and Struve functions}

\begin{theorem}
Let $\nu>0$, $z\in\C\setminus(-\infty,0]$. Then
\begin{eqnarray}
\int_0^\infty J_\nu(cz) J_{2\nu}\left(2\sqrt{c}\right) dc
= \frac{1}{z} J_\nu\left(\frac{1}{z}\right).
\end{eqnarray}
\end{theorem}
\medskip

\noindent {\bf Proof.}\quad We are given the integral \cite[(6.514.1)]{Grad}
\[
\int_{0}^{\infty} J_{\nu}\left( \frac{a}{x} \right) J_{\nu}(bx) dx 
= b^{-1} J_{2\nu} \left( 2\sqrt{ab} \right),
\]
where ${\rm Re}\, \nu>0$, $a$,$b>0$.
By applying Theorem \ref{2:Hankel} to the function
$F_\nu^a:(0,\infty)\to\C$ defined by 
$
F_\nu^a(x):=x^{-1} J_{\nu} \left( a x^{-1} \right)
$,
we obtain the following integral
\[
\int_{0}^{\infty}
  J_{\nu}(bx) J_{2\nu} \left( 2\sqrt{ab} \right)
  db
= x^{-1} J_{\nu} \left( \frac{a}{x} \right),
\]
where $\nu>0$, $a>0$, $x\in\C\setminus(-\infty,0]$. 
By making the substitutions $x=az$, $c=ba$, we obtain the desired result.
$\hfill\blacksquare$

\begin{theorem}
Let $-\frac12 \leq \nu < \frac52$, $z\in\C\setminus(-\infty,0]$. Then
\begin{eqnarray}
\int_{0}^{\infty} c J_\nu(cz)
&&\left[e^{i (\nu+1)\pi/2} K_{2\nu} \left(2e^{i \pi/4} \sqrt{c}\right)
+e^{-i (\nu+1)\pi/2} K_{2\nu} \left(2e^{-i \pi/4} \sqrt{c}\right)\right]
dc
\nonumber\\
&&= \frac{1}{z^3} K_\nu\left(\frac{1}{z}\right).
\end{eqnarray}
\end{theorem}
\medskip

\noindent {\bf Proof.}\quad We are given the integral \cite[(6.514.3)]{Grad}
\[
\int_{0}^{\infty} J_\nu \left(\frac{a}{x}\right) K_\nu(bx)dx
= b^{-1} e^{i (\nu+1)\pi/2} K_{2\nu} \left[2e^{i \pi/4} \sqrt{ab}\right]
+b^{-1} e^{-i (\nu+1)\pi/2} K_{2\nu} \left[2e^{-i \pi/4} \sqrt{ab}\right],
\]
where $a>0$, ${\rm Re}\, b>0$, $|{\rm Re}\, \nu|<\frac{5}{2}$.
By applying Theorem \ref{2:Hankel} to the function
$F_\nu^b:(0,\infty)\to\C$ defined by
$
F_\nu^b(x):= bx^{-3} K_\nu(bx^{-1}),
$
we obtain the following integral
\[
\int_{0}^{\infty} a J_\nu(ax)
\left[e^{i (\nu+1)\pi/2} K_{2\nu} \left(2e^{i \pi/4} \sqrt{ab}\right)
+e^{-i (\nu+1)\pi/2} K_{2\nu} \left(2e^{-i \pi/4} \sqrt{ab}\right)\right]
da
= \frac{b}{x^3} K_\nu\left(\frac{b}{x}\right),
\]
where ${\rm Re}\, b>0$, $-\frac12 \leq \nu < \frac52$, $x\in\C\setminus(-\infty,0]$.
With the substitutions $x=bz$, $c=ba$, we obtain the desired result.
$\hfill\blacksquare$

\begin{theorem}
Let $|\nu|<\frac{1}{2}$, $z\in\C\setminus(-\infty,0]$. Then
\begin{eqnarray}
\int_{0}^{\infty} J_\nu(cz) 
\left[K_{2\nu}\left(2\sqrt{c}\right) 
- \frac{\pi}{2}Y_{2\nu}\left(2\sqrt{c}\right)\right]dc
= -\frac{\pi}{2z} Y_\nu\left(\frac{1}{z}\right).
\end{eqnarray}
\end{theorem}
\medskip

\noindent {\bf Proof.}\quad We are given the integral \cite[(6.514.4)]{Grad}
\[
\int_{0}^{\infty} Y_\nu\left(\frac{a}{x}\right) J_\nu(bx) dx
= -\frac{2b^{-1}}{\pi} \left[K_{2\nu}\left(2\sqrt{ab}\right) 
- \frac{\pi}{2}Y_{2\nu}\left(2\sqrt{ab}\right)\right],
\]
where $a$,$b>0$, $|{\rm Re}\, \nu|<\frac{1}{2}$.
By applying Theorem \ref{2:Hankel} to the function
$F_\nu^a:(0,\infty)\to\C$ defined by
\[
F_\nu^a(x):= -\frac{\pi}{2x} Y_\nu\left(\frac{a}{x}\right),
\]
we obtain the following integral
\[
\int_{0}^{\infty} J_\nu(bx) 
\left[K_{2\nu}\left(2\sqrt{ab}\right) 
- \frac{\pi}{2}Y_{2\nu}\left(2\sqrt{ab}\right)\right]db
= -\frac{\pi}{2x} Y_\nu\left(\frac{a}{x}\right),
\]
where $a>0$, $|\nu|<\frac{1}{2}$, $x\in\C\setminus(-\infty,0]$.
With the substitutions $x=az$, $c=ab$, we obtain the desired result.
$\hfill\blacksquare$

\begin{theorem}
Let $\nu\geq-\frac14$, $\mu>-\frac12$, $z\in\C\setminus(-\infty,0]$. Then
\begin{eqnarray}
\int_{0}^{\infty} J_{2\nu} (cz) J_{\nu} \left( \frac{c^2}{4} \right) c\,dc 
= 2 J_{\nu}\left(z^2\right),
\end{eqnarray}
\begin{eqnarray}
\int_{0}^{\infty} J_{\mu}(cz) J_{\mu} \left( \frac{1}{4c} \right) dc
= z^{-1} J_{2\mu} \left(\sqrt{z}\right).
\end{eqnarray}
\end{theorem}
\medskip

\noindent {\bf Proof.}\quad We are given the integral \cite[(6.516.1)]{Grad}
\[
\int_{0}^{\infty} J_{2\nu} \left(a\sqrt{x}\right) J_{\nu}(bx) dx 
= b^{-1} J_{\nu} \left( \frac{a^2}{4b} \right),
\]
where ${\rm Re}\, \nu>-\frac12$, $a$,$b>0$.
By applying Theorem \ref{2:Hankel} to the functions
$F_\nu^b:(0,\infty)\to\C$ and $G_\mu^a:(0,\infty)\to\C$ defined by 
$
F_\nu^b(x):=2 b J_{\nu}\left(bx^2\right)
$,
$
G_\mu^a(x):=x^{-1} J_{2\mu} \left(a\sqrt{x}\right)
$,
we obtain the following integrals
\[
\int_{0}^{\infty} a J_{2\nu} (ax) J_{\nu} \left( \frac{a^2}{4\beta} \right) da 
= 2 \beta J_{\nu}\left(\beta x^2\right),
\]
\[
\int_{0}^{\infty} J_\mu(bx) J_\mu\left(\frac{\alpha^2}{4b}\right) db
= x^{-1} J_{2\mu}(\alpha\sqrt{x}),
\]
where $\alpha$,$\beta>0$, $\nu\geq-\frac14$, $\mu>-\frac12$, $z\in\C\setminus(-\infty,0]$. 
With the substitutions $z^2=bx^2$, $c=a/\sqrt{b}$, and $\sqrt{z}=a\sqrt{x}$, $c=b/a^2$
respectively, we obtain the desired results.
$\hfill\blacksquare$
\medskip

\begin{theorem}
Let $\nu>-1$, $\mu\geq-\frac12$, $z\in\C\setminus(-\infty,0]$. Then
\begin{eqnarray}
\int_0^\infty J_{\nu/2}(cz) J_{\nu/2}\left(\frac{1}{4c}\right) dc
=z^{-1} J_\nu\left(\sqrt{z}\right),
\end{eqnarray}
\begin{eqnarray}
\int_0^\infty c J_\mu(cz) J_{\mu/2}\left(\frac{c^2}{4}\right) dc
=2J_{\mu/2}(z^2).
\end{eqnarray}
\end{theorem}
\medskip

\noindent {\bf Proof.}\quad We are given the integral \cite[(6.526.1)]{Grad}
\[
\int_0^\infty x J_{\nu/2}(ax^2) J_\nu(bx) dx
=\frac{1}{2a} J_{\nu/2}\left(\frac{b^2}{4a}\right),
\]
where $a,b>0$, ${\rm Re}\,\nu>-1$. 
By applying Theorem \ref{2:Hankel} to the function
$F_\nu^b:(0,\infty)\to\C$ defined by 
$
F_\nu^b(x):=x^{-1} J_\nu\left(b\sqrt{x}\right)
$,
we obtain the following integrals
\[
\int_0^\infty J_{\nu/2}(ax) J_{\nu/2}\left(\frac{\beta^2}{4a}\right) da
=x^{-1} J_\nu\left(\beta\sqrt{x}\right),
\]
\[
\int_0^\infty b J_\mu(bx) J_{\mu/2}\left(\frac{b^2}{4\alpha}\right) db
=2\alpha J_{\mu/2}(\alpha x^2),
\]
where $\alpha$,$\beta>0$, $\nu>-1$, $\mu\geq-\frac12$, $x\in\C\setminus(-\infty,0]$.
With the substitutions $\sqrt{z}=\beta\sqrt{x}$, $c=a/\beta^2$, and $z^2=\alpha x^2$, $x=z\sqrt{\alpha}$, 
we obtain the desired results.
$\hfill\blacksquare$
\medskip

\begin{theorem}
Let $\nu\geq-\frac14$, $x\in\C\setminus(-\infty,0]$. Then
\begin{eqnarray}
\int_0^\infty a^2 J_{2\nu}(ax) J_{\nu+1/2}(a^2) da
= \frac x4 J_{\nu-1/2}\left(\frac{x^2}4\right).
\end{eqnarray}
\end{theorem}
\medskip

\noindent {\bf Proof.}\quad By applying Theorem \ref{2:Hankel} to the function
$F_\nu:(0,\infty)\to\C$ defined by
$
F_\nu(x):=\frac x4 J_{\nu-1/2}\left(\frac{x^2}4\right)
$,
we obtain the desired result from the known integral \cite[(6.527.1)]{Grad}
\[
\int_0^\infty J_{2\nu}(2ax) J_{\nu-1/2}(x^2) dx
=\frac12 a \,J_{\nu+1/2}(a^2),
\]
where $a>0$, ${\rm Re}\,\nu>-\frac12$.
$\hfill\blacksquare$

\begin{theorem}
Let $\nu\geq-\frac14$, $x\in\C\setminus(-\infty,0]$. Then
\begin{eqnarray}
\int_0^\infty a^2 J_{2\nu}(ax) J_{\nu-1/2}(a^2) da
= \frac x4 J_{\nu+1/2}\left(\frac{x^2}4\right).
\end{eqnarray}
\end{theorem}
\medskip

\noindent {\bf Proof.}\quad By applying Theorem \ref{2:Hankel} to the function
$F_\nu:(0,\infty)\to\C$ defined by
$
F_\nu(x):=\frac x4 J_{\nu+1/2}\left(\frac{x^2}4\right)
$,
we obtain the desired result from the known integral \cite[(6.527.1)]{Grad}
\[
\int_0^\infty J_{2\nu}(2ax) J_{\nu+1/2}(x^2) dx
=\frac12 a \,J_{\nu-1/2}(a^2),
\]
where $a>0$, ${\rm Re}\,\nu>-2$.
$\hfill\blacksquare$
\medskip

\begin{theorem}
Let $\nu\geq-\frac12$, $z\in\C\setminus(-\infty,0]$. Then
\begin{eqnarray}
\int_0^\infty c J_\nu(cz) \emph{\textbf{H}}_{\nu/2}\left(\frac{c^2}{4}\right) dc
= -2 Y_{\nu/2}(z^2).
\end{eqnarray}
\end{theorem}
\medskip

\noindent {\bf Proof.}\quad We are given the integral \cite[(6.526.4)]{Grad}
\[
\int_0^\infty x Y_{\nu/2}(ax^2) J_\nu(bx) dx
= -\frac{1}{2a} \textbf{H}_{\nu/2}\left(\frac{b^2}{4a}\right),
\]
where $a>0$, ${\rm Re}\,b>0$, ${\rm Re}\,\nu>-1$ and $\textbf{H}_\nu:\C\to\C$, 
for $\nu\in\N_0$, is 
the Struve function defined in \cite[(11.2.1)]{NIST}.
By applying Theorem \ref{2:Hankel} to the function
$F_\nu^a:(0,\infty)\to\C$ defined by 
$
F_\nu^a(x):=-2a Y_{\nu/2}(ax^2)
$,
we obtain the following integral
\[
\int_0^\infty b J_\nu(bx) \textbf{H}_{\nu/2}\left(\frac{b^2}{4a}\right) db
= -2a Y_{\nu/2}(ax^2),
\]
where $a>0$, $\nu\geq-\frac12$, $x\in\C\setminus(-\infty,0]$.
With the substitutions $z^2=ax^2$, $c=b/\sqrt{a}$, we obtain the desired result.
$\hfill\blacksquare$
\medskip

\section{Exponential, logarithmic and inverse trigonometric functions}

\begin{theorem}
Let $\nu\geq-\frac12$, $z\in\C\setminus(-\infty,0]$. Then
\begin{eqnarray}
\int_0^\infty 
J_\nu(cz) e^{-2/c} c^{-1} dc
= 2 J_\nu(2\sqrt{z}) K_\nu(2\sqrt{z}).
\end{eqnarray}
\end{theorem}
\medskip

\noindent {\bf Proof.}\quad We are given the integral \cite[(6.526.4)]{Grad}
\[
\int_0^\infty x J_\nu(2\sqrt{ax}) K_\nu(2\sqrt{ax}) J_\nu(bx) dx
=\frac12 b^{-2} e^{-2a/b},
\]
where ${\rm Re}\,a>0$, $b>0$, ${\rm Re}\,\nu>-1$.
By applying Theorem \ref{2:Hankel} to the function
$F_\nu^a:(0,\infty)\to\C$ defined by
$
F_\nu^a(x):=2 J_\nu(2\sqrt{ax}) K_\nu(2\sqrt{ax})
$,
we obtain the following integral
\[
\int_0^\infty b^{-1} J_\nu(bx) e^{-2a/b} db
= 2 J_\nu(2\sqrt{ax}) K_\nu(2\sqrt{ax}),
\]
where ${\rm Re}\,a>0$, $\nu\geq-\frac12$, $x\in\C\setminus(-\infty,0]$.
With the substitutions $z=ax$, $c=b/a$, we obtain the desired result.
$\hfill\blacksquare$
\medskip

\begin{theorem}
Let $a>0$, $z\in\C\setminus(-\infty,0]$. Then
\begin{eqnarray}
\int_{0}^{a} J_{1}(bz) \ln \left( 1-\frac{b^2}{a^2} \right) db
= -\pi z^{-1} Y_0(az).
\end{eqnarray}
\end{theorem}
\medskip

\noindent {\bf Proof.}\quad We apply Theorem \ref{2:Hankel} to the function
$F_a:(0,\infty)\to\C$ defined by
$
F_a(x):= -\pi x^{-1} Y_0(ax)
$,
where $Y_\nu:\C\setminus(-\infty,0]\to\C$, (order) $\nu\in\C$, is the Bessel 
function of the second kind defined in \cite[(10.2.3)]{NIST}. We obtain the 
desired result from the known integral \cite[(6.512.6)]{Grad}
\[
\int_{0}^{\infty} J_1(bx) Y_0(ax) dx = -\frac{b^{-1}}{\pi} \ln \left(1-
\frac{b^2}{a^2} \right),
\]
where $0<b<a$.
$\hfill\blacksquare$

\begin{theorem}
Let $z\in\C\setminus(-\infty,0]$. Then
\begin{eqnarray}
\int_0^\infty J_1(cz) \ln(1+c^2) dc
= 2 z^{-1} K_0(z).
\end{eqnarray}
\end{theorem}
\medskip

\noindent {\bf Proof.}\quad We are given the integral \cite[(6.512.9)]{Grad}
\[
\int_{0}^{\infty} K_0(ax) J_1(bx) dx = \frac{1}{2b} \ln \left(1+\frac{b^2}
{a^2}\right),
\]
where $a$,$b>0$ and $K_\nu:\C\setminus(-\infty,0]\to\C$, (order) $\nu\in\C$, 
is the modified Bessel function of the second kind defined in \cite[(10.25.3)]{NIST}. 
We apply Theorem \ref{2:Hankel} to the function
$F_a:(0,\infty)\to\C$ defined by
$
F_a(x):=2x^{-1}K_0(ax)
$,
and obtain the following integral
\[
\int_0^\infty J_1(bx) \ln \left(1+\frac{b^2}{a^2}\right) db
= 2 x^{-1} K_0(ax),
\]
where $a>0$, $x\in\C\setminus(-\infty,0]$.
With the substitutions $z=ax$ and $c=b/a$, we obtain the desired result.
$\hfill\blacksquare$

\begin{theorem}
Let $a>0$, $z\in\C\setminus(-\infty,0]$. Then
\begin{eqnarray}
\int_0^{2a} \sin^{-1} \left(\frac{b}{2a}\right) J_1(bz) db
= \frac{\pi}{2z} [J_0^2(az)-J_0(2az)].
\end{eqnarray}
\end{theorem}
\medskip

\noindent {\bf Proof.}\quad By applying Theorem \ref{2:Hankel} to the function
$F_a:(0,\infty)\to\C$ defined by 
\[
F_a(x):=\frac{\pi}{2x} J_0^2(ax),
\]
we obtain the desired result from the known integral \cite[(6.513.9)]{Grad}
\[
\int_{0}^{\infty} J_0^2(ax) J_1(bx)dx=\left\{
\begin{array}{ll}
    {\displaystyle b^{-1}} & $if$\,\, 0<2a<b, \\[.1cm]
    {\displaystyle \frac{2}{\pi b}\sin^{-1} \left(\frac{b}{2a}\right)} & $if$\,\, 0<b<2a.
\end{array}
  \right.
\]
$\hfill\blacksquare$

\section{Hypergeometric and Legendre functions}

\begin{theorem}
Let 
$n\in\N_0$, 
$\mu>0$, $\nu>\frac12$, 
$t\in\C\setminus (-\infty,0]$. Then
\begin{eqnarray}
\int_{0}^{\alpha} 
J_{\nu-n-1}(bt) \,{}_2F_1\left(\begin{array}{c} \nu,-n\\
\nu-n\end{array}; \frac{b^2}{\alpha^2} \right) 
b^{\nu-n}\,
db
=
  \frac{n!\alpha^{\nu-n}\Gamma(\nu-n)J_{\nu+n}(\alpha t)}{t\Gamma(\nu)},
\end{eqnarray}
\begin{eqnarray}
\int_{\beta}^{\infty} 
J_{\mu+n}(a t) \,{}_2F_1\left(\begin{array}{c} \mu,-n\\
\mu-n\end{array}; \frac{\beta^2}{a^2} \right) 
a^{-\mu+n+1} \,
da
=
  \frac{n!\beta^{-\mu+n+1}\Gamma(\mu-n)J_{\mu-n-1}(\beta t)}{t\Gamma(\mu)}.
\end{eqnarray}
\end{theorem}
\medskip

\noindent {\bf Proof.}\quad By applying Theorem \ref{2:Hankel} to the functions
$G_n^{\nu,\alpha}:(0,\infty)\to\C$ and $H_n^{\mu,\beta}:(0,\infty)\to\C$ defined by 
\[
G_n^{\nu,\alpha}(t):=\frac{n!\alpha^{\nu-n}\Gamma(\nu-n)J_{\nu+n}(\alpha t)}{t\Gamma(\nu)},
\]
\[
H_n^{\mu,\beta}(t):=\frac{n!\beta^{-\mu+n+1}\Gamma(\mu-n)J_{\mu-n-1}(\beta t)}{t\Gamma(\mu)},
\]
we obtain the desired results from the known integral \cite[(6.512.2)]{Grad}
\[
\int_{0}^{\infty} J_{\nu+n}(\alpha t) J_{\nu-n-1}(\beta t) dt = \left\{
\begin{array}{ll}
    \frac{\beta^{\nu-n-1} \Gamma(\nu)}{\alpha^{\nu-n}n!\,\Gamma(\nu-n)} 
    \,{}_2F_1\left(\begin{array}{c}\nu,-n\\
    \nu-n\end{array}; \frac{\beta^2}{\alpha^2} \right)
    & $if$\,\, 0<\beta<\alpha, \\[.2cm]
    (-1)^n (2\alpha)^{-1} & $if$\,\, \beta=\alpha, \\[.2cm]
    0 & $if$\,\, \beta>\alpha,
\end{array}
\right.
\]
where ${\rm Re}\,\nu>0$ and ${}_2F_1:\C^2\times(\C\setminus-\N_0)\times(\C\setminus[1,\infty))\to\C$ 
is the hypergeometric function 
defined 
in \cite[(15.2.1)]{NIST}. 
$\hfill\blacksquare$

\begin{theorem}
Let $\nu\geq-\frac12$, $\nu>-2\,{\rm Re}\,\mu-1$, $z\in\C\setminus(-\infty,0]$. Then
\begin{eqnarray}
\int_{0}^{\infty} 
P_{-1/2+\nu/2}^{-\mu}\left(\sqrt{1+\frac{4}{c^2}}\right) 
&&Q_{-1/2+\nu/2}^{-\mu}\left(\sqrt{1+\frac{4}{c^2}}\right) 
J_\nu(cz) dc
\nonumber\\
&&= \frac{e^{-\mu\pi i}\,\Gamma\left(\frac{\nu-2\mu+1}{2}\right)}{z
\Gamma\left(\frac{\nu+2\mu+1}{2}\right)} I_\mu(z) K_\mu(z).
\end{eqnarray}
\end{theorem}
\medskip

\noindent {\bf Proof.}\quad We are given the integral \cite[(6.513.3)]{Grad}
\[
\int_{0}^{\infty} I_\mu(ax) K_\mu(ax) J_\nu(bx) dx
= \frac{e^{\mu\pi i} \Gamma\left(\frac{\nu+2\mu+1}{2}\right)}
{b \Gamma\left(\frac{\nu-2\mu+1}{2}\right)} 
P_{-1/2+\nu/2}^{-\mu}\left(\sqrt{1+\frac{4a^2}{b^2}}\right) 
Q_{-1/2+\nu/2}^{-\mu}\left(\sqrt{1+\frac{4a^2}{b^2}}\right),
\]
where ${\rm Re}\, a>0$, $b>0$, ${\rm Re}\, \nu>-1$, ${\rm Re}\, \nu+2\mu>-1$ and
$P_\nu^\mu:\C$ \textbackslash $(-\infty,1]\to\C$, for $\nu+\mu\notin-\N$ with degree $\nu$ and 
order $\mu$, and $Q_\nu^\mu:\C$ \textbackslash $(-\infty,1]\to\C$, for $\nu+\mu\notin-\N$ with 
degree $\nu$ and order $\mu$, are the associated Legendre functions of the
first \cite[(14.3.6)]{NIST} and second \cite[(14.3.7)]{NIST} kind respectively. 
By applying Theorem \ref{2:Hankel} to the function
$F_\nu^{\mu,a}:(0,\infty)\to\C$ defined by
\[
F_\nu^{\mu,a}(x):= \frac{e^{-\mu\pi i}\,\Gamma\left(\frac{\nu-2\mu+1}{2}\right)}{x
\Gamma\left(\frac{\nu+2\mu+1}{2}\right)} I_\mu(ax) K_\mu(ax),
\]
we obtain the following integral
\[
\int_{0}^{\infty} 
P_{-1/2+\nu/2}^{-\mu}\left(\sqrt{1+\frac{4a^2}{b^2}}\right) 
Q_{-1/2+\nu/2}^{-\mu}\left(\sqrt{1+\frac{4a^2}{b^2}}\right) 
J_\nu(bx) db
= \frac{e^{-\mu\pi i}\,\Gamma\left(\frac{\nu-2\mu+1}{2}\right)}{x
\Gamma\left(\frac{\nu+2\mu+1}{2}\right)} I_\mu(ax) K_\mu(ax),
\]
where ${\rm Re}\, a>0$, $\nu\geq-\frac12$, $\nu>-{\rm Re}\,2\mu-1$, $x\in\C\setminus(-\infty,0]$.
By making the substitutions $z=ax$, $c=b/a$, we obtain the desired result.
$\hfill\blacksquare$

\begin{theorem}
Let $\nu>\mp{\rm Re}\,\mu-1$, $\nu\geq-\frac12$, $z\in\C\setminus(-\infty,0]$. Then
\begin{eqnarray}
\int_{0}^{\infty} J_\nu(cz) \left[Q_{-1/2+\nu/2}^{-\mu}
\left(\sqrt{1+\frac{4}{c^2}}\right)\right]^2 dc
= \frac{e^{-2\mu\pi i} \Gamma\left(\frac{1+\nu-2\mu}{2}\right)}
{z \Gamma\left(\frac{1+\nu+2\mu}{2}\right)} [K_\mu(z)]^2.
\end{eqnarray}
\end{theorem}
\medskip

\noindent {\bf Proof.}\quad We are given the integral \cite[(6.513.5)]{Grad}
\[
\int_{0}^{\infty} [K_\mu(ax)]^2 J_\nu(bx) dx
= \frac{e^{2\mu\pi i} \Gamma\left(\frac{1+\nu+2\mu}{2}\right)}
{b \Gamma\left(\frac{1+\nu-2\mu}{2}\right)} 
\left[Q_{-1/2+\nu/2}^{-\mu}
\left(\sqrt{1+\frac{4a^2}{b^2}}\right)\right]^2,
\]
where ${\rm Re}\, a>0$, $b>0$, ${\rm Re}\,(\frac\nu2 \pm \mu)>-\frac12$.
By applying Theorem \ref{2:Hankel} to the function
$F_\nu^{\mu,a}:(0,\infty)\to\C$ defined by
\[
F_\nu^{\mu,a}(x):= \frac{e^{-2\mu\pi i} \Gamma\left(\frac{1+\nu-2\mu}{2}\right)}
{x \Gamma\left(\frac{1+\nu+2\mu}{2}\right)} [K_\mu(ax)]^2,
\]
we obtain the following integral
\[
\int_{0}^{\infty} J_\nu(bx) \left[Q_{-1/2+\nu/2}^{-\mu}
\left(\sqrt{1+\frac{4a^2}{b^2}}\right)\right]^2 db
= \frac{e^{-2\mu\pi i} \Gamma\left(\frac{1+\nu-2\mu}{2}\right)}
{x \Gamma\left(\frac{1+\nu+2\mu}{2}\right)} [K_\mu(ax)]^2,
\]
where ${\rm Re}\, a>0$, $\nu>\mp{\rm Re}\,\mu-1$, $\nu\geq-\frac12$, $x\in\C\setminus(-\infty,0]$.
With the substitutions $z=ax$, $c=b/a$, we obtain the desired result.
$\hfill\blacksquare$

\section{Jacobi polynomials and Chebyshev polynomials of the first kind}

\begin{theorem}
Let $n\in\N_0$, $\nu>-n-1$, $\alpha$,$\beta>0$, 
$z\in\C\setminus(-\infty,0]$. Then
\begin{eqnarray}
\int_{\beta}^{\infty} 
P_{n}^{(\nu,0)}\left(1-\frac{2\beta^2}{a^2}\right) J_{\nu+2n+1}(az)
a^{-\nu} 
\,da
= z^{-1} \beta^{-\nu} J_\nu(\beta z),
\end{eqnarray}
\begin{eqnarray}
\int_{0}^{\alpha} 
P_{n}^{(\nu,0)} \left(1-\frac{2b^2}{\alpha^2}\right) J_\nu(bz)
b^{\nu+1}
\, db
= z^{-1} \alpha^{\nu+1} J_{\nu+2n+1}(\alpha z).
\end{eqnarray}
\end{theorem}
\medskip

\noindent {\bf Proof.}\quad By applying Theorem \ref{2:Hankel} to the functions
$F_\nu^b:(0,\infty)\to\C$ and $G_\nu^{a,n}:(0,\infty)\to\C$ defined by 
$
F_\nu^b(x):=x^{-1} b^{-\nu} J_\nu(bx)
$,
$
G_\nu^{a,n}(x):=x^{-1} a^{\nu+1} J_{\nu+2n+1}(ax)
$,
we obtain the desired results from the known integral \cite[(6.512.4)]{Grad}
\[
\int_{0}^{\infty} J_{\nu+2n+1}(ax) J_\nu(bx) dx=\left\{
\begin{array}{ll}
    b^\nu a^{-\nu-1} P_{n}^{(\nu,0)} \left(1-2a^{-2}b^2\right)
    & $if$\,\, 0<b<a, \\[.2cm]
    0 & $if$\,\, 0<a<b,
\end{array}
\right.
\]
where ${\rm Re}\,\nu>-n-1$ and $P_{n}^{(\alpha,\beta)}:\C\to\C$ is the Jacobi polynomial 
defined in \cite[(18.3.1)]{NIST}.
$\hfill\blacksquare$

\begin{theorem}
Let $a>0$, $\nu\geq-\frac12$, $z\in\C\setminus(-\infty,0]$. Then
\begin{eqnarray}
\int_0^{2a} \frac{J_\nu(bz)}{\sqrt{4a^2-b^2}} T_n\left(\frac{b}{2a}\right) db
= \frac{\pi}{2} J_{(\nu+n)/2}(az) J_{(\nu-n)/2}(az). 
\end{eqnarray}
\end{theorem}
\medskip

\noindent {\bf Proof.}\quad By applying Theorem \ref{2:Hankel} to the function
$F_\nu^{n,a}:(0,\infty)\to\C$ defined by 
\[
F_\nu^{n,a}(x)=\frac{\pi}{2} J_{(\nu+n)/2}(ax) J_{(\nu-n)/2}(ax),
\]
we obtain the desired result from the known integral \cite[(6.522.11)]{Grad}
\[
\int_0^\infty xJ_{(\nu+n)/2}(ax) J_{(\nu-n)/2}(ax) J_\nu(bx) dx
=
\left\{
\begin{array}{ll}
    2\pi^{-1}b^{-1}(4a^2-b^2)^{-1/2}\, T_n(\frac{b}{2a}) & $if$\,\, 0<b<2a,\\[.2cm]
    0 & $if$\,\, 2a<b, 
\end{array}
\right.
\]
where ${\rm Re}\,\nu>-1$ and $T_n:\C\to\C$, for $n\in\N_0$, is the Chebyshev polynomial of the first kind 
found in \cite[(18.3.1)]{NIST}.
$\hfill\blacksquare$
\medskip

\section{Examples where the Hankel transforms fails}
\label{AttemptFailure}

In the following
examples, the potential use of the method given by Theorem \ref{2:Hankel}
fails because the condition (\ref{2:cond}) can not be satisfied.
This has been verified by analyzing the well-understood behavior of the integrands 
in a small neighborhood of the endpoints $\{0,\infty\}$.
\begin{itemize}
\item The definite integral \cite[(6.512.1)]{Grad} with
$
G_\nu^{\mu,b}(x):=
\alpha(\nu,\mu)\Gamma(\nu+1)b^{-\nu}x^{-1}J_\nu(bx),
$
$
H_\nu^{\mu,a}(x):=
\alpha(\nu,\mu)
\Gamma(\nu+1)a^{\nu+1}x^{-1}J_\mu(ax),
$
where
$\alpha(\nu,\mu):=\Gamma\left(\frac{\mu-\nu+1}{2}\right)
/\Gamma\left(\frac{\mu+\nu+1}{2}\right)$.\\[-0.6cm]

\item 
The definite integral \cite[(6.514.1)]{Grad} 
with 
$
G_\nu^b(x):=b x^{-3} J_{\nu} ( bx^{-1} ).
$\\[-0.6cm]

\item The definite integral \cite[(6.514.2)]{Grad} with 
$
F_\nu^b(x):=b x^{-3} Y_{\nu}( bx^{-1}).
$\\[-0.6cm]

\item The definite integrals \cite[(6.516.2)]{Grad}, \cite[(6.516.3)]{Grad}, 
\cite[(6.516.4)]{Grad}, and \cite[(6.516.7)]{Grad} with respectively
$
F_\nu^b(x):=-2b Y_{\nu}(bx^2)
$,
$
F_\nu^b(x):=4b\pi^{-1} K_\nu(bx^2)
$,
$
F_\nu^a(x):=x^{-1} Y_{2\nu}(a\sqrt{x})
$,
and
$
F_\nu^a(x):=\frac{4 \pi^{-1} x^{-1}}{\sec(\nu \pi)} K_{2\nu}(a \sqrt{x}).
$\\[-0.6cm]

\item The definite integral \cite[(6.522.2)]{Grad} with
$
F_\nu^{\mu,a}(x):=\tfrac12 e^{-2\mu\pi i}
\beta(\nu,\mu) [K_\mu(ax)]^2,
$
where
$\beta(\nu,\mu):=\Gamma(\frac\nu2-\mu)/\Gamma(1+\frac\nu2+\mu).$\\[-0.6cm]

\item The definite integrals \cite[(6.522.6)]{Grad} and \cite[(6.522.8)]{Grad} 
with 
$F_a(x):=-\frac{\pi}{2} J_0(ax) Y_0(ax),$
and
$
G_\nu^{\mu,a}(x):=
\tfrac12 e^{-2\mu\pi i}
\beta(\nu,\mu)
K_{\mu-1/2}(ax) K_{\mu+1/2}(ax),
$ respectively.\\[-0.6cm]

\item The definite integral \cite[(6.522.16)]{Grad} with
$
F_\nu^{b,c}(x):= \sqrt{\pi} 
x^\nu\gamma(\nu)
I_\nu(cx) K_\nu(bx)
$, where
$\gamma(\nu):=(8bc)^{-\nu}/\Gamma\left(\nu+\frac{1}{2}\right).$\\[-0.6cm]

\item The definite integrals \cite[(6.526.2)]{Grad} and \cite[(6.526.3)]{Grad} 
with
$F_\nu^b(x):=2x^{-1} Y_\nu(b\sqrt{x}),$
$G_\nu^b(x):= \cos(\frac{\nu \pi}{2})K_\nu(b\sqrt{x})/(2\pi x),$
respectively.\\[-0.6cm]

\item The definite integral \cite[(6.526.6)]{Grad} with
$
F_\nu^a(x):=4a \pi^{-1} K_{\nu/2}(ax^2).
$\\[-0.6cm]

\item The definite integral \cite[(6.527.3)]{Grad} with
$
F_\nu(x):=-x Y_{\nu+1/2}\left(x^2/4\right)/4.
$\\[-0.6cm]
\end{itemize}




\begin{thebibliography}{1}

\bibitem{CohlVolkmerDefInt}
H.~S. Cohl and H.~{Volkmer}.
\newblock Definite integrals using orthogonality and integral transforms.
\newblock {\em Symmetry, Integrability and Geometry:~Methods and Applications},
  8, 2012.

\bibitem{NIST:DLMF}
{NIST Digital Library of Mathematical Functions}.
\newblock {Release 1.0.9 of 2014-08-29}.
\newblock {Online companion to \cite{NIST}}.

\bibitem{NIST}
{F.~W.~J. Olver and D.~W. Lozier and R.~F. Boisvert and C.~W. Clark}, editor.
\newblock {\em {{NIST} {H}andbook of {M}athematical {F}unctions}}.
\newblock {Cambridge University Press}, {New York, NY}, {2010}.
\newblock {Print companion to \cite{NIST:DLMF}}.

\bibitem{Grad}
I.~S. Gradshteyn and I.~M. Ryzhik.
\newblock {\em Table of {I}ntegrals, {S}eries, and {P}roducts}.
\newblock Elsevier/Academic Press, Amsterdam, seventh edition, 2007.

\bibitem{Ostermannetal12}
A.~Ostermann and G.~Wanner.
\newblock {\em Geometry by its {H}istory}.
\newblock Undergraduate Texts in Mathematics. Readings in Mathematics.
  Springer, Heidelberg, 2012.

\bibitem{Stempak}
K.~Stempak.
\newblock A new proof of {S}onine's formula.
\newblock {\em Proceedings of the American Mathematical Society},
  104(2):453--457, 1988.

\bibitem{Watson}
G.~N. Watson.
\newblock {\em A {T}reatise on the {T}heory of {B}essel {F}unctions}.
\newblock Cambridge Mathematical Library. Cambridge University Press,
  Cambridge, second edition, 1944.

\end{thebibliography}

\def\cprime{$'$} \def\dbar{\leavevmode\hbox to 0pt{\hskip.2ex \accent"16\hss}d}

\end{document}